\documentclass[11pt]{article}
\usepackage[greek,english]{babel}
\usepackage[iso-8859-7]{inputenc}
\usepackage{amssymb}
\usepackage{amsmath}
\usepackage{amsthm}
\usepackage{latexsym}
\usepackage{amsfonts}
\usepackage{graphicx}
\usepackage{graphics}

\theoremstyle{plain}
\newtheorem{thm}{Theorem}[section]   
\newtheorem{prop}[thm]{Proposition}

\theoremstyle{definition}

\newtheorem*{Proof}{Proof}

\newcommand{\ld} {{\ldots}}

\newcommand{\thi} {{\theta}}
\newcommand{\Thi} {{\varTheta}}
\newcommand{\de} {{\delta}}

\newcommand{\si} {{\sigma}}

\newcommand{\la} {{\lambda}}

\newcommand{\e} {{\varepsilon}}

\newcommand{\dis}{\displaystyle}

\newcommand{\ct}{{\cal{T}}}

\newcommand{\ch}{{\cal{H}}}

\newcommand{\tn}{{\widetilde{n}}}

\newcommand{\ra}{{\rightarrow}}

\newcommand{\qb}{$\quad\blacksquare$}
\def\1{\it1\hspace*{-0.150cm}{\footnotesize{I}}}

\def\C{{\mathbb{C}}}
\def\Q{{\mathbb{Q}}}
\def\N{{\mathbb{N}}}

\begin{document}
\title{\bf Simultaneous approximation of\\ translation operators}
\author{\bf N. Tsirivas}
%
\date{}
\maketitle
\noindent
\begin{abstract}
Let $(m_n)_{n\in\N}$ be an unbounded sequence of complex numbers and $(\thi_v)_{v\in\N}$ be a sequence of distinct numbers in $[0,1)$ where $\N$ is the set of natural numbers. For every $v, n\in\N$ we shall consider the translation functions $t_{v,n}:\C\ra\C$, where $t_{v,n}(z)=z+m_ne^{2\pi i\thi_v}$ for every $v,n\in\N$, $z\in\C$, where $\C$ is the set of complex numbers.

We shall also consider the linear and continuous translation, operators:
\[
T_{v,n}:\ch(\C)\ra\ch(\C), \ \ \text{so that}
\]
\[
T_{v,n}(f)=f\circ t_{v,n} \ \ \text{for every} \ \ v,n\in\N, \ \ f\in\ch(\C).
\]
Let $\rho$ be the usual metric in $\ch(\C)$ that defines the topology of uniform convergence on compacta.

We fix an entire function $g$.

We shall prove that there is a pair $(x,\la_n)$, where $x\in\ch(\C)$ and $(\la_n)$ is a sequence of natural numbers so that the pair $(x,\la_n)$ is a solution of system of equations:
\[
\lim_{n\ra+\infty}\rho(T_{v,m_{\la_n}}(x),g)=0, \ \ \text{for all} \ \ v\in\N.
\]
\end{abstract}
 \smallskip\\
MSC (2020) 47A16

\section{Introduction}\label{sec1}
\noindent

A classical result of Birkhoff \cite{2}, which goes back to 1929, says that there are entire functions of which the integer translates are dense in the space of all entire functions, endowed with the topology $\ct_u$ of local uniform convergences (see also Luh \cite{8} for a more general statement). Birkhoff's proof was constructive.

Much later, during the 80's, Gethner and Shapiro \cite{6} and independently Grosse-Erdmann \cite{7} showed that Birkhoff's result can be recovered as a particular case of a much more general theorem through the use of Baire's category theorem. This approach simplified Birkhoff's argument substantially and, in addition it gave us precise information on the topological size of these functions. In particular, Grosse-Erdmann proved that for every fixed sequence of complex numbers $(w_n)$ with $w_n\ra\infty$, the set $\{f\in\ch(\C)|\overline{\{f(z+w_n):n\in\N\}}=\ch(\C)\}$ is $G_\de$ and dense in $\ch(\C)$, and hence topologicallty ``large''.

Let us apply this result, in certain cases.

Let $(\thi_v)_{v\in\N}$ be a sequence of distinct numbers in $[0,1)$ and $(m_n)_{n\in\N}$ be a sequence of complex numbers so that $m_n\ra\infty$. We shall consider the numbers $w_n(\thi_v)=m_ne^{2\pi i\thi_v}$, $n,v\in\N$. That is, for every $v\in\N$ we shall consider the sequence $(w_n(\thi_v))_{n\in\N}$. Of course, we have $w_n(\thi_v)\ra\infty$ as $n\ra+\infty$ for every $v\in\N$.

We now set:
\[
E_v=\{f\in\ch(\C)|\overline{\{f(z+w_n(\thi_v)):n\in\N\}}=\ch(\C)\} \ \ \text{for every} \ \ v\in\N.
\]
Based on Grosse-Erdmann's result we conclude that for every $v\in\N$ the set $E_v$ is $G_\de$ and dense in $\ch(\C)$.

So, the set $E=\dis\bigcap^{+\infty}_{v=1}E_v$ is a $G_\de$ dense subset of $\ch(\C)$, so it is non-empty by Baire's Category Theorem, given that the space $\ch(\C)$ is a complete metric space. Let us see in more detail what this result means.

Let $f\in E$. Then, for every $v\in\N$ and $g\in\ch(\C)$ there is a subsequence $\la^v_n=\la_n(v,g)$ of $(w_n(\thi_v))$, that depends on $g$ and $v$ so that for every compact set $K\subseteq\C$
\[
\sup_{z\in K}|f(z+\la^v_n)-g(z)|\ra0 \ \ \text{as} \ \ n\ra+\infty.
\]
So, this convergence depends on the specific sequence $\la^v_n=\la_n(v,g)$, $n\in\N$, and the sequence $\la^v_n$ depends on the specific  number $\thi_v\in[0,1)$. In the present paper we shall examine whether we can have this convergence without the dependence on the specific number $\thi_v\in[0,1)$.

Thus, we shall introduce the set of entire functions that achieve simultaneous approximation on all numbers $\thi_v$, $v\in\N$, where $\thi_v\in[0,1)$ for every $v\in\N$ with the same sequence of indices.

Given that, we shall consider the set:\\
$SA=($\,by the initial letters of the words Simultaneous Approximations)\\
$\hspace*{0.65cm}=\{f\in\ch(\C)|$ for every $g\in\ch(\C)$ there is a sequence $(\la_n)_{n\in\N}$, so that \\
$\la_n\in\{m_n,n\in\N\}$ for every $n\in\N$ and, subsequently, for every $K\subseteq\C$, $K$ compact and $v\in\N$\\
$\dis\sup_{z\in K}|f(z+\la_ne^{2\pi i\thi_v})-g(z)|\ra0$ as $n\ra+\infty\}$.\\
Of course $SA\subseteq E$.

We prove that the set $SA$ is a $G_\de$-dense subset of $\ch(\C)$, so it is non-empty. In order to prove that $SA$ is a $G_\de$, dense subset of $\ch(\C)$ we introduce one other set $V\subseteq\ch(\C)$ and we prove that $V$ is a $G_\de$, dense subset of $\ch(\C)$ and $SA=V$.

There are some papers concerning common hypercyclic vectors for translation operators see \cite{10}, \cite{11}, \cite{1}, \cite{9}, \cite{5}, \cite{3}, \cite{4}.

Whenever, we refer to a topology in the $\ch(\C)$ space, we mean the topology of uniform convergence on compacta.

In the following Section \ref{sec2} we prove some helpful propositions in order to prove our main result Theorem \ref{thm2.6}.
\section{The main result}\label{sec2}
\noindent

First of all, we shall prove a proposition which is the key in order to prove our main result.

We fix $g\in\ch(\C)$.

We fix some natural numbers $n_0\ge2$, $v_0,N_0$, and some real numbers $\thi_1,\thi_2,\ld,\thi_{n_0}$ where $\thi_i\in[0,1)$ for each $i=1,\ld,n_0$ and $\thi_i\neq\thi_j$ for every $i,j\in T_{n_0}=\{1,\ld,n_0\}$, $i\neq j$.

For every natural number $m$ we use the set
\[
V_g(m,v_0,N_0,n_0)\!=\!\Big\{f\in\ch(\C)\Big/\sup_{|z|\le v_0}\Big|f(z+me^{2\pi i\thi_j})-g(z)\Big|\!<\!\frac{1}{N_0} \; \text{for every} \;  j\!=\!1,\ld,n_0\Big\}.
\]
For every $m\in\N$, $j\in T_{n_0}$ we use the set
\[
V_g(m,v_0,N_0,j)=\Big\{f\in\ch(\C)\Big|\sup_{|z|\le v_0}\Big|f(z+me^{2\pi\thi_j})-g(z)\Big|<\frac{1}{N_0}\Big\}.
\]
Of course, we have
\begin{eqnarray}
V_g(m,v_0,N_0,n_0)=\bigcap^{n_0}_{j=1}V_g(m,v_0,N_0,j),  \label{eq1}
\end{eqnarray}
based on the above definitions.

It is easy to see that the sets $V_g(m,v_0,N_0,j)$ are open in $\ch(\C)$ for every $m\in\N$, $j=1,\ld,n_0$, so the set $V_g(m,v_0,N_0,n_0)$ is open in $\ch(\C)$ for every $m\in\N$, according to the above relation (\ref{eq1}).

Therefore, the set $\dis\bigcup^{+\infty}_{m=1}V_g(m,v_0,N_0,n_0)$ is open in $\ch(\C)$.

We shall prove now the following proposition.
\begin{prop}\label{prop2.1}
According to the above notations the set $\dis\bigcup^{+\infty}_{m=1}V_g(m,v_0N_0,n_0)$ is dense in $\ch(\C)$.
\end{prop}
\begin{Proof}
We fix that $h\in\ch(\C)$, a compact set $K\subseteq\C$ and $\e>0$. It suffices to conclude that there is $f\in\ch(\C)$ and $m_0\in\N$, so that
\setcounter{equation}{0}
\begin{eqnarray}
f\in V_g(m_0,v_0,N_0,n_0) \ \  \text{and} \ \  \|f-h\|_K<\si  \label{eq1}
\end{eqnarray}
We set $D_v=\{z\in\C\mid|z|\le v\}$ for every $v\in\N$. We also choose $v_1\in\N$ so that
\begin{eqnarray}
D_{v_0}\cup K\subseteq D_{v_1}.  \label{eq2}
\end{eqnarray}
Let's assume that $m\in\N$, so that
\[
D_{v_1}\cap(D_{v_1}+me^{2\pi i\thi_j})\neq\emptyset
\]
for some $j\in T_{n_0}$ (if it exists).\\
This means that there also exist $z_j,w_j\in D_{v_1}$, so that
\begin{eqnarray}
w_j=z_j+me^{2\pi i\thi_j} \ \ \text{for some} \ \ j\in T_{n_0}.  \label{eq3}
\end{eqnarray}
According to (\ref{eq3}), we shall have:
\[
|w_j-z_j|=m \ \ \text{and this gives} \ \ m\le 2v_1.
\]
Therefore, for every $m\in\N$ and $m>2v_1$, we have
\begin{eqnarray}
D_{v_1}\cap(D_{v_1}+me^{2\pi i\thi_j})=\emptyset \ \ \text{for every} \ \ j\in T_{n_0}.  \label{eq4}
\end{eqnarray}
Let $j_1,j_2\in T_{n_0}$, so that $j_1\neq j_2$. (We remind $T_{n_0}=\{1,2,\ld,n_0\}$).

Let $m\in\N$ so that
\[
(D_{v_1}+me^{2\pi i\thi_{j_1}})\cap(D_{v_1}+me^{2\pi\thi_{j_2}})\neq\emptyset \ \ \text{(if it exists)}.
\]
This means that there are $z_1,w_1\in D_{v_1}$ so that
\begin{eqnarray}
z_1+me^{2\pi i\thi_{j_1}}=w_1+me^{2\pi i\thi_{j_2}}. \label{eq5}
\end{eqnarray}
By (\ref{eq5}) we have:
\begin{eqnarray}
|z_1-w_1|=m\big|e^{2\pi i(\thi_{j_2}-\thi_{j_1})}-1\big|.  \label{eq6}
\end{eqnarray}
By (\ref{eq6}) we have:
\begin{eqnarray}
m\le\frac{2v_1}{\big|e^{2\pi i(\thi_{j_2}-\thi_{j_1})}-1\big|}.  \label{eq7}
\end{eqnarray}
So, for every $m\in\N$, so that
\[
m>\frac{2v_1}{\big|e^{2\pi i(\thi_{j_2}-\thi_{j_1})}-1\big|} \ \ \text{we have}:
\]
\begin{eqnarray}
(D_{v_1}+me^{2\pi i\thi_{j_1}})\cap(D_{v_1}+me^{2\pi i\thi_{j_2}})=\emptyset.  \label{eq8}
\end{eqnarray}
We set
\[
M_0=\min\big\{\big|e^{2\pi i(\thi_{j_2}-\thi_{j_1})}-1\big|,j_1,j_2\in T_{n_0},j_1\neq j_2\big\}.
\]
We fix now some natural number $m_0$ so that $m_0>\max\Big\{2v_1,\dfrac{2v_1}{M_0}\Big\}$. Then, by (\ref{eq4}) and (\ref{eq8}) we shall take that:
\[
D_{v_1}\cap(D_{v_1}+m_0e^{2\pi i\thi_j})=\emptyset \ \ \text{for every} \ \ j\in T_{n_0}
\]
and
\begin{eqnarray}
(D_{v_1}+m_0e^{2\pi i\thi_{j_1}})\cap(D_{v_1}+m_0e^{2\pi i\thi j_2})=\emptyset \ \  \text{for every} \ \ j_1,j_2\in T_{n_0}, \ \ j_1\neq j_2.  \label{eq9}
\end{eqnarray}
We set:
\[
L=D_{v_1}\cup\bigg(\bigcup^{n_0}_{j=1}(D_{v_1}+m_0e^{2\pi i\thi_j})\bigg).
\]
Because of (\ref{eq9}) we have that the set $L$ is a union of $n_0+1$ disjoint closed discs with the same radius $v_1$.

This means that the set $L$ is a compact set with connected complement. We shall consider the function $F:L\ra\C$, defined as follows:
\[
F(z)=\left\{\begin{array}{l}
              h(z), \; \text{given that} \; z\in D_{v_1} \\ [1ex]
              g(z-m_0e^{2\pi i\thi_j}) \\ [1ex]
              \text{if}\; z\in D_{v_1}+m_0e^{2\pi i\thi_j} \;\text{for some}\;j\in T_{n_0}.
            \end{array}\right.
\]
Of course, $F$ is continuous on  $L$ and holomorphic on $\overset{\circ}{L}$. So, according to Mergelyan's Theorem there is a complex polynomial $f$, so that
\begin{eqnarray}
\|F-f\|_L<\min\bigg\{\e,\frac{1}{N_0}\bigg\}.  \label{eq10}
\end{eqnarray}
Based on the definition of $F$ and (\ref{eq10}), we have:
\begin{eqnarray}
\|f-h\|_K<\e,  \label{eq11}
\end{eqnarray}
because of relation (\ref{eq2}), and the definition of $L$.

Let's suppose $w\in D_{v_1}$. Then, for $j\in T_{n_0}$, $w+m_0e^{2\pi i\thi_j}\in D_{v_1}+m_0e^{2\pi i\thi_j}$.

We also set $z=w+m_0e^{2\pi i\thi_j}\in D_{v_1}+m_0e^{2\pi i\thi_j}$, $j\in T_{n_0}$.

Then, $F(z)=g(z-m_0e^{2\pi i\thi_j})=g(w)$. By (\ref{eq10}) we have for every $w\in D_{v_1}$ and $j\in T_{n_0}$,
\[
\big|f(w+m_0e^{2\pi i\thi_j})-g(w)\big|<\frac{1}{N_0}.
\]
This means that we shall have:
\[
\big|f(w+m_0e^{2\pi i\thi_j})-g(w)\big|<\frac{1}{N_0} \ \ \text{for every} \ \ w\in D_{v_1} \ \ \text{and} \ \ j\in T_{n_0}.
\]
This concludes that $f\in V_g(m_0,v_0,N_0,n_0)$, because $f\in\ch(\C)$ as a polynomial and relation (\ref{eq2}). By this fact and (\ref{eq11}) the proof of this proposition is complete now because relation (\ref{eq1}) is satisfied. \qb
\end{Proof}
Now, we shall fix a sequence $(m_s)_{s\in\N}$ of complex numbers, that is unbounded. With the notation of the previous Proposition \ref{prop2.1} we shall consider the set:
\[
V_g(m_s,v_0,N_0,n_0)=\bigg\{f\in\ch(\C)\Big|\sup_{|z|\le v_0}\Big|f(z+m_se^{2\pi i\thi_j})-g(z)\bigg|<\frac{1}{N_0} \ \ \text{for every} \ \ j\in T_{n_0}\bigg\},
\]
for every $s\in\N$.

The sets $V_g(m_s,v_0,N_0,n_0)$ are open for every $s\in\N$, so the set $\dis\bigcup^{+\infty}_{s=1}V_g(m_s,v_0,N_0,n_0)$ is open in $\ch(\C)$.

As in Proposition \ref{prop2.1}, we shall now prove the following proposition:
\begin{prop}\label{prop2.2}
The set $\dis\bigcup^{+\infty}_{s=1}V_g(m_s,v_0,N_0,n_0)$ is dense in $\ch(\C)$.
\end{prop}
\begin{Proof}
The proof is similar to that of Proposition \ref{prop2.1} and for this reason the proof is omitted. \qb
\end{Proof}

The space $\ch(\C)$ is separable.

Let $(p_k)_{k\in\N}$ be a fixed dense sequence of $\ch(\C)$ (for example $(p_k)_{k\in\N}$ be an enumeration of all complex polynomials with coefficients in $\Q+i\Q$). For every $v,N,k,n,s\in\N$, $n\ge2$ we shall consider the set:
\[
V_{p_k}(m_s,v,N,n)=\bigg\{f\in\ch(\C)\bigg|\sup_{|z|\le v}\Big|f(z+m_s
e^{2\pi i\thi_j})-p_k(z)\bigg|<\frac{1}{N} \ \ \text{for every} \ \ j\in T_n\bigg\}
\]
The sets $V_{p_k}(m_s,v,N,n)$ are open in $\ch(\C)$ for every $v,N,k,n,s\in\N$, $n\ge2$, so that the set
$\dis\bigcup^{+\infty}_{s=1}V_{p_k}(m_s,v,N,n)$ is open for every $v,N,k,n\in\N$, $n\ge2$. According to Proposition \ref{prop2.2}, we shall have that the sets $\dis\bigcup^{+\infty}_{s=1}V_{p_k}(m_s,v,N,n)$ are dense in $\ch(\C)$ for every $v,N,k,n\in\N$, $n\ge2$.

We shall also consider the set:
\[
V=\bigcap^{+\infty}_{v=1}\bigcap^{+\infty}_{N=1}\bigcap^{+\infty}_{k=1}\bigcap^{+\infty}_{n=2}
\bigg(\bigcup^{+\infty}_{s=1}V_{p_k}(m_s,v,N,n)\bigg).
\]
By the above notation, we shall examine the following proposition:
\begin{prop}\label{prop2.3}
The set $V$ is a $G_\de$-dense subset of $\ch(\C)$, so $V$ is non-empty.
\end{prop}
\begin{Proof}
The set $V$ is a $G_\de$ subset of $\ch(\C)$, by its definition, because the sets\linebreak $\dis\bigcup^{+\infty}_{s=1}V_{p_k}(m_s,v,N,n)$ are open for every $v,N,k,n\in\N$, $n\ge2$. Based on Proposition \ref{prop2.2}, the sets $\dis\bigcup^{+\infty}_{s=1}V_{p_k}(m_s,v,N,n)$ are dense for every $v,N,k,n\in\N$, $n\ge2$, so the result holds by Baire's Category Theorem because the space $\ch(\C)$ is a complete metric\linebreak space. \qb
\end{Proof}

We connect now the previous set $V$ with the set of entire functions that succeed simultaneous approximation in a countable set of real numbers.

We shall state here the respective data. Let $(\thi_n)_{n\in\N}$ be a sequence of real numbers, so that $\thi_n\in[0,1)$ and $\thi_{j_1}\neq\thi_{j_2}$ for every $j_1,j_2\in\N$, $j_1\neq j_2$, $n\in\N$. Let $(m_s)_{s\in\N}$ be a fixed sequence of complex numbers which is unbounded.

Let $\Thi=\{x\in[0,1)|$ there is $n\in\N$ so that $x=\thi_n\}$. Of course, the set $\Thi$ is infinite, and the set:
\[
m=\{y\in\C\,|, \text{there is some}\; s\in\N\; \text{such that}\; y=m_s\}
\]
is infinite also.

We shall consider the set:
\[
SA=\{f\in\ch(\C)\,|\; \text{for every} \; g\in\ch(\C),
\]
there is a sequence $(\la_n)_{n\in\N}$ so that $\la_n\in m$ for every $n\in\N$, so that for every $a\in\Thi$ and for every compact set $K\subseteq\C$
\[
\sup_{z\in K}\big|f(z+\la_ne^{2\pi ia})-g(z)\big|\ra0 \ \ \text{as} \ \ n\ra+\infty\}.
\]
Our aim is to prove that the set $SA$ is non empty.

The method to prove that $SA\neq\emptyset$ is the following:\\
We shall prove that $SA=V$ and given that $V\neq\emptyset$ we shall also have $SA\neq\emptyset$.

In order to prove that $SA=V$ we prove that $SA\subseteq V$ and $V\subseteq SA$. This is the subject of the following two propositions.
\begin{prop}\label{prop2.4}
It holds $SA\subseteq V$.
\end{prop}
\begin{Proof}
If $SA=\emptyset$, then the result is proven. We suppose that $SA\neq\emptyset$. Let $f\in SA$. We fix $v_0,N_0,n_0,k_0\in\N$, $n_0\ge2$.

Because $f\in SA$ for $g=p_{k_0}$ there is a sequence $(\la_n)_{n\in\N}$, so that $\la_n\in m$ for every $n\in\N$, so that for every $a\in\Thi$ and compact set $K\subseteq\C$
\[
\sup_{z\in K}|f(z+\la_ne^{2\pi ia})-p_{k_0}(z)|\ra0 \ \ \text{as} \  n\ra+\infty.
\]
So, for $K=D_{v_0}$ we have that
\[
\sup_{|z|\le v_0}|f(z+\la_ne^{2\pi i\thi_j})-p_{k_0}(z)|\ra0 \ \ \text{as} \ \ n\ra+\infty,
\]
for every $j\in T_{n_0}$.

This concludes that for every $j\in T_{n_0}$ there is some $n_j\in\N$, so that
\[
\sup_{|z|\le v_0}\big|f(z+\la_ne^{2\pi i\thi})-p_{k_0}(z)\Big|<\frac{1}{N_0} \ \ \text{for every} \ \ n\in\N, \ \ n\ge n_j.
\]
Let $\tn=\max\{n_j|j\in T_{n_0}\}$.\\
We shall take that
\[
\sup_{|z|\le v_0}\big|f(z+\la_ne^{2\pi i\thi_j})-p_{n_0}(z)\Big|<\frac{1}{N_0}
\]
for every $j\in T_{n_0}$, for every $n\in\N$, $n\ge\tn$.

This concludes that $f\in V_{p_{k_0}}(\la_\tn,v_0,N_0,n_0)$, or else $f\in\dis\bigcup^{+\infty}_{s=1}V_{p_{k_0}}(m_s,v_0,N_0.n_0)$ because $\la_\tn\in m$, that implies $f\in V$ and the result is proven. \qb
\end{Proof}
\begin{prop}\label{prop2.5}
It holds that $V\subseteq SA$.
\end{prop}
\begin{Proof}
We know that $V\neq\emptyset$. Let $f\in V$. We shall prove that $f\in SA$.

We fix $g\in\ch(\C)$. We shall show that there exists a sequence $(\la_n)_{n\in\N}$, so that $\la_n\in m$ for every $n\in\N$, so that for every $a\in\Thi$ and compact set $K\subseteq\C$
\[
\sup_{z\in K}\big|f(z+\la_ne^{2\pi ia})-g(z)\big|\ra0 \ \ \text{as} \ \ n\ra+\infty.
\]
Based on the above mentioned, we shall construct now the respective sequence $(\la_n)_{n\in\N}$. We shall fix some $n_0\in\N$, $n_0\ge2$.

Given that the sequence $(p_k)_{k\in\N}$ of complex polynomials with coefficients in $\Q+i\Q$ is dense in $\ch(\C)$, there is some $k_0\in\N$ so that
\setcounter{equation}{0}

\begin{eqnarray}
\|g-p_{k_0}\|_{D_{n_0}}<\frac{1}{2n_0}.  \label{eq1}
\end{eqnarray}
Because $f\in V$ we shall have: $f\in\dis\bigcup^{+\infty}_{s=1}V_{p_{k_0}}(m_s,n_0,2n_0,n_0)$. This means that there is some $s_{n_0}\in\N$ so that $f\in V_{p_{k_0}}(m_{s_{n_0}},n_0,2n_0,n_0)$, or else
\begin{eqnarray}
\sup_{|z|\le n_0}\big|f(z+m_{s_{n_0}}e^{2\pi i\thi_j})-p_{k_0}(z)\Big|<\frac{1}{2n_0}, \ \ \text{for every} \ \ j\in T_{n_0}.  \label{eq2}
\end{eqnarray}
By (\ref{eq1}), (\ref{eq2}) and the triangle inequality we have:
\begin{eqnarray}
\sup_{|z|\le n_0}\big|f(z+m_{s_{n_0}}e^{2\pi i\thi_j})-g(z)\Big|<\frac{1}{n_0} \ \ \text{for every} \ \ j\in T_{n_0}. \label{eq3}
\end{eqnarray}
According to the previous procedure for every $n\in\N$, $n\ge2$, we shall choose some $s_n\in\N$, so that:
\begin{eqnarray}
\sup_{|z|\le n}\big|f(z+m_{s_n}e^{2\pi i\thi_j})-g(z)\Big|<\frac{1}{n} \ \ \text{for every} \ \ j\in T_n. \label{eq4}
\end{eqnarray}
We shall prove now that for the sequence $(s_n)_{n\in\N}$ has as follows:\\
$\dis\sup_{z\in K}\Big|f(z+m_{s_n}e^{2\pi ia})-g(z)\Big|\ra0$ as $n\ra+\infty$, for every compact set $K\subseteq\C$ and for every $a\in\Thi$.

We fix some $\e_0>0$.

There are $v_0\in\N$ and $n_0\in\N$, so that $K\subseteq D_v$ for every $v\in\N$, $v\ge v_0$ and $a_0=\thi_{n_0}$. Let suppose $N_0\in\N$, so that $\dfrac{1}{N_0}<\e_0$.

Let $M_0=\max\{v_0,n_0,N_0,2\}$. For every $n\in\N$, $n\ge M$, we have $n\ge v_0$, so $K\subseteq D_n$. Of course, $a_0\in\{\thi_1,\thi_2,\ld,\thi_n\}$ for every $n\in\N$, $n\ge M_0$, because $a_0=\thi_{n_0}$ and $n_0\le M_0\le n$. We shall also have:
\[
\frac{1}{n}\le\frac{1}{M_0}\le\frac{1}{N_0}<\e_0 \ \ \text{for every} \ \ n\in\N, \ \ n\ge M_0.
\]
So, for every $n\in\N$, $n\ge M_0$ we shall have:
\[
\sup_{z\in K}\Big|f(z+m_{s_n}e^{2\pi ia_0})-g(z)|\le\sup_{|z|\le n}\big|
f(z+m_{s_n}e^{2\pi ia_0})-g(z)\Big|<\frac{1}{n}<\e_0
\]
by (\ref{eq4}). This concludes that
\[
\sup_{z\in K}\big|f(z+m_{s_n}e^{2\pi ia_0})-g(z)\big|\ra0 \ \ \text{as} \ \ n\ra+\infty.
\]
So, for every $a\in\Thi$ and compact set $K\subseteq\C$ we have:
\[
\sup_{z\in K}\big|f(z+m_{s_n}e^{2\pi ia})-g(z)\big|\ra0 \ \ \text{as} \ \ n\ra+\infty.
\]
Because this is the case for arbitrary $g\in\ch(\C)$ we conclude that $f\in SA$ and the proof of this proposition is complete. \qb
\end{Proof}

Based on the above results, we are ready now to state and prove the main result of this paper, that is Theorem \ref{thm2.6}.
\begin{thm}\label{thm2.6}
The set $SA$ is a $G_\de$ dense subset of $\ch(\C)$. More specifically, the set $SA$ is non-empty.
\end{thm}
\begin{Proof}
Based on Proposition \ref{prop2.4} and \ref{prop2.5} we have that $SA=V$. We have also proved in Proposition \ref{prop2.3} that the set $V$ is a $G_\de$ and dense subset of $\ch(\C)$. So, the result follows. \qb
\end{Proof}
\vspace*{1cm}
N. Tsirivas  \\
University of Thessaly, \\
Department of Mathematics, \\
Lamia, Greece.\\
email: ntsirivas@uth.gr


\begin{thebibliography}{xxx}
\bibitem{1} F. Bayart, Common hypercyclic vectors for high-dimensional families of operators. Int. Math. Res. Notices 21 (2016), 6512-6552.
%
\bibitem{2} G. D. Birkhoff, D\'{e}monstration, d' un th\'{e}or\`{e}m \'{e}t\`{e}mentaire sur les fonctions enti\`{e}res, C. R. Acad. Sci. Paris, 189 (1929), 473-475.
%
\bibitem{3} G. Costakis, Approximation by translates of entire functions, in Complex and Harmonic Analysis, Destech, Publ. Inc.
%
\bibitem{4} G. Costakis, M. Sambarino, Genericity of wild holomorphic functions and common hypercyclic vectors. Adv. Math. 182 (2004), 278, 306.
%
\bibitem{5} G. Costakis, N. Tsirivas, V. Vlachou, Non-existence of common hypercyclic vectors for certain families of translation operators, Comput. Methods Funct. Theory, 15 (2015), 393-401.
%
\bibitem{6} R. M. Gethner, J. H. Shapiro, Universal vectors for operators on spaces of holomorphic functions, Proc. Amer, Math. Soc. 100 (1987), 281-288.
%
\bibitem{7} K. G. Grosse-Erdmann, Holomorphe Monster und universelle Funktionen [Holomorphic monsters and universal functions] [German], Ph.D. Dissertation, University of Trier, Trier 1987, Mitt. Math. Sem. Giessen, vol. 176, 1987.
%
\bibitem{8} W. Luh, On universal functions, Colloq Math. Soc. J\`{a}nos Bolyai 19 (1976), 503-511.
%
\bibitem{9} S. Shkarin, Remarks on common hypercyclic vectors, J. Funct. Anal. 258 (2010), 132-160.
%
\bibitem{10} N. Tsirivas, Common hypercyclic functions for translation operators with large gaps, J. Funct. Anal. 272 (2017), 2726-2751.
%
\bibitem{11} N. Tsirivas, Existence of common hypercyclic vectors for translation operators, Journal of Operator Theory, 80, 2 (2018), 257-294.
\end{thebibliography}
\end{document}